\documentclass[preprint,12pt,review]{elsarticle}

\usepackage{amssymb,latexsym}
\usepackage{amsmath}

\biboptions{sort&compress}

\journal{: \;  J. Number Theory}

\begin{document}

\newtheorem{lema}{Lemma}
\newtheorem{teo}{Theorem}
\newproof{prova}{Proof}

\begin{frontmatter}

\title{A shortened recurrence relation for the Bernoulli numbers}

\author{F. M. S. Lima}

\address{Institute of Physics, University of Brasilia, P.O. Box 04455, 70919-970, Brasilia-DF, Brazil}


\ead{fabio@fis.unb.br}

\date{\today}

\begin{abstract}
In this note, starting with a little-known result of Kuo, I derive a recurrence relation for the Bernoulli numbers $B_{2 n}$, $n$ being any positive integer. This new recurrence seems advantageous in comparison to other known formulae since it allows the computation of both $B_{4 n}$ and $B_{4 n +2}$ from only $B_0, B_2, \ldots, B_{2n}$.
\newline
\end{abstract}


\end{frontmatter}


\vspace{0.5cm}

The Bernoulli numbers $B_n$, $n$ being a nonnegative integer, can be defined by the generating function
\begin{equation}
\frac{x}{\mathrm{e}^x -1} = \sum_{n=0}^\infty{B_n \, \frac{x^n}{n!}} \, , \quad |x| < 2 \pi.
\label{eq:DefBn}
\end{equation}
The first few values are well-known: $B_0 = 1$, $B_1 = -\frac12$, and $B_2 = \frac16$. It is also well-known that $B_n = 0$ for odd values of $n$, $n>1$. For even values of $n$ ($=2m$), $n>1$, the numbers $B_{2m}$ form a subsequence of non-null real numbers such that $(-1)^{m+1} \, B_{2m}>0$. In other words, the entire sequence of Bernoulli numbers up to any $B_n$, $n > 1$, consists of $B_0$ and $B_1$, given above, and the preceding numbers $B_{2m}$, $m = 1, \ldots, \lfloor n/2 \rfloor$.  The basic properties of these numbers can be found in Sec.~9.61 of~\cite{Gradsh}.

The numbers $B_n$ appear in many instances in pure and applied mathematics, most notably in number theory, finite differences calculations, and asymptotic analysis. Therefore, the efficient computation of the numbers $B_n$ is of great interest. To this end, recurrence formulae were soon recognized as the most efficient tool~\cite{livro_old}. One of the simplest such relations is found by multiplying both sides of Eq.~\eqref{eq:DefBn} by $\mathrm{e}^x-1$, using the Cauchy product with the Maclaurin series for $\mathrm{e}^x-1$, and equating the coefficients of the powers of $x$, which results in
\begin{equation}
\sum_{j=0}^n{\binom{n+1}{j} \, B_j} = 0 \, , \quad n \ge 1 \, .
\label{eq:Bn_simples}
\end{equation}
This kind of recurrence formula has the disadvantage of demanding the previous knowledge of all $B_0, B_1, \ldots, B_{n-1}$ for the computation of $B_n$. In searching for more efficient formulae, shortened recurrence relations of two different types have been discovered. The first type consists of the so-called lacunary recurrence relations, in which $B_n$ is determined only from every second, or every third, etc., preceding Bernoulli numbers (see, e.g., the lacunary formula by Ramanujan~\cite{Ramanujan}).  The second type demands the knowledge of only the second-half of the Bernoulli numbers up to $B_{n-1}$ in order to compute $B_n$~\cite{Agoh}. For an extensive study of these and other recurrence relations involving the Bernoulli numbers, see~\cite{livrao}.

Here in this note, I apply the Euler's formula relating the even zeta value $\zeta(2n)$ to $B_{2n}$ to a little-known recurrence formula for $\zeta{(2n)}$ obtained by Kuo~\cite{Kuo}, in order to convert it into a recurrence formula for $B_{2n}$. By doing this, in fact I introduce a third type of recurrence relation for the Bernoulli numbers, in the sense that it allows us to compute both $B_{4 n}$ and $B_{4 n +2}$ from only the first-half of the preceding numbers, i.e. $B_0, B_2, \ldots, B_{2n}$. Then, the efficiency of this new recurrence formula certainly surpasses that of most known formulae, mainly for large values of $n$.

For real values of $s$, $\,s>1$, the Riemann zeta function is defined as $\zeta(s) := \sum_{n=1}^\infty{{\,1/n^s}}$.  In this domain, the convergence of this series is guaranteed by the integral test.\footnote{For $s=1$, one has the harmonic series $\sum_{n=1}^\infty{{1/n}}$, which diverges to infinity.}  For positive even values of $s$, one has the well-known Euler's formula (1740)~\cite{Euler}:
\begin{equation}
\zeta{(2 n)} = (-1)^{n-1} \, \frac{2^{2n-1} \, B_{2n}}{(2n)!} \: \pi^{2n} \, ,
\label{eq:Euler}
\end{equation}
which, in face of the rationality of every Bernoulli number $B_n$, yields $\zeta{(2 n)}$ as a rational multiple of $\pi^{2 n}$.  The function $\zeta{(s)}$, as defined above, can be extended to the entire complex plane (except the only simple pole at $s=1$) by analytic continuation, which yields $\zeta{(0)} = -\frac12$.\footnote{Alternatively, we can use the globally convergent series $\frac{1}{1-2^{1-s}} \, \sum_{n=0}^\infty {\frac{1}{2^{n+1}} \, \sum_{k=0}^n {(-1)^k \, \binom{n}{k}} \, \frac{1}{(k+1)^s}}$, due to Hasse (1930)~\cite{Hasse}, which is valid for all complex numbers $s \ne 1 + \frac{2 \pi m}{\ln{2}} \, \mathrm{i}$, $m$ being any integer. At $s=0$, it reduces to $\zeta{(0)} = -\sum_{n=0}^\infty {\frac{1}{2^{n + 1}} \, \sum_{k=0}^n {(-1)^k \, \binom{n}{k}} } = -\frac12 -\frac12 \sum_{n=1}^\infty {\frac{1}{2^n} \sum_{k=0}^n {(-1)^k \, \binom{n}{k}}}$. Of course, $\sum_{k=0}^n {(-1)^k \, \binom{n}{k}} = 0$ for all $n>0$, so $\zeta{(0)} = -\frac12$.}  The reader should note that Eq.~\eqref{eq:Euler} remains valid for $n=0$.

These are the necessary ingredients to state our first lemma, which comes from a little-known recurrence formula by Kuo (1949)~\cite{Kuo}, just written directly in terms of the Riemann zeta function.

\vspace{0.5cm}

\begin{lema} [Kuo's recurrence formula for $\zeta{(2n)}$]
\label{lem:Kuo}
\quad For any positive integer $n$, one has
\begin{eqnarray*}
\zeta{(2n)} = \frac{2^{2n-1 \, \pi^{2n}}}{4 \, {(n-1)!}^2 \, (2n-1)} +\frac{1}{(n-1)!} \, \sum_{k=0}^{\lfloor n/2 \rfloor}{(-1)^k \, \frac{\zeta{(2k)}\,(2 \pi)^{2n-2k}}{(n-2k)! \, (2n-2k)}} \nonumber \\
+\frac{1}{\pi} \, \sum_{k=0}^{\lfloor n/2 \rfloor}{\sum_{j=0}^{\lfloor n/2 \rfloor}{(-1)^{k+j} \, \zeta{(2k)} \, \zeta{(2j)} \, \frac{(2 \pi)^{2n-2k-2j+1}}{(n-2k)! \, (n-2j)! \, (2n-2k-2j+1)} } } \, .
\label{eq:Kuo}
\end{eqnarray*}
\end{lema}

\vspace{0.5cm}

This recurrence formula is proved in~\cite{Kuo} by developing successive integrations (from $0$ to $x$) of the Fourier series $\sum_{n=1}^\infty{\sin{(n\,x)}/n} = {\,(\pi -x)/2}$, which converges for all positive real $x < 2\,\pi$, and then applying the Parseval's theorem.

We are now in a position to prove the following theorem.

\vspace{0.5cm}

\begin{teo}[Recurrence relation for $B_{2n}$]
\label{teo:B2n}
\quad For any positive integer $n$, one has
\begin{eqnarray}
B_{2n} = (-1)^{n-1} \, \Biggl[ a_n -b_n \, \sum_{k=0}^{\lfloor n/2 \rfloor}{\frac{B_{2k}}{(2k)! \, (n-2k)! \, (n-k)}} \nonumber \\
 +(2n)! \, \sum_{k=0}^{\lfloor n/2 \rfloor}{\frac{B_{2k}}{(2k)! \, (n-2k)!} \, \sum_{j=0}^{\lfloor n/2 \rfloor}{\frac{B_{2j}}{(2j)! \, (n-2j)!} \: \frac{1}{2n-2k-2j+1} } } \Biggl] \, ,
\end{eqnarray}
where $a_n = \frac{n}{2} \, \frac{(2n-2)!}{(n-1)!^2}$ and $b_n = \frac{(2n)!}{2\,(n-1)!}$.
\end{teo}

\begin{prova}
\; By dividing both sides of the Kuo's recurrence formula, as given in Lemma~\ref{lem:Kuo}, by $\pi^{2n}$, one has
\begin{eqnarray}
\frac{\zeta{(2n)}}{\pi^{2n}} = \frac{2^{2n-1 \, \pi^{2n}}}{4 \, {(n-1)!}^2 \, (2n-1)} +\frac{1}{(n-1)!} \, \sum_{k=0}^{\lfloor n/2 \rfloor}{(-1)^k \, \frac{\zeta{(2k)}}{\pi^{2k}} \, \frac{2^{2n-2k}}{(n-2k)! \, (2n-2k)}} \nonumber \\
+\sum_{k=0}^{\lfloor n/2 \rfloor}{\sum_{j=0}^{\lfloor n/2 \rfloor}{(-1)^{k+j} \, \frac{\zeta{(2k)}}{\pi^{2k}} \, \frac{\zeta{(2j)}}{\pi^{2j}} \, \frac{2^{2n-2k-2j+1}}{(n-2k)! \, (n-2j)! \, (2n-2k-2j+1)} } } \, .
\label{eq:Kuo2}
\end{eqnarray}
From Euler's equation, Eq.~\eqref{eq:Euler}, one knows that $\frac{\zeta{(2m)}}{\pi^{2m}} = (-1)^{m-1} \, \frac{2^{2m-1} \, B_{2m}}{(2m)!}$, which is valid for all integer $m$, $m \ge 0$. By substituting this in Eq.~\eqref{eq:Kuo2}, one finds, after some algebra,
\begin{eqnarray}
\frac{B_{2n}}{(2n)!} = \frac{(-1)^{n-1}}{4 \, {(n-1)!}^2 \, (2n-1)} -\frac{(-1)^{n-1}}{(n-1)!} \, \sum_{k=0}^{\lfloor n/2 \rfloor}{\frac{B_{2k}}{(2k)!} \, \frac{1}{(n-2k)! \, (2n-2k)}} \nonumber \\
+ (-1)^{n-1} \, \sum_{k=0}^{\lfloor n/2 \rfloor}{\sum_{j=0}^{\lfloor n/2 \rfloor}{\frac{B_{2k}}{(2k)!} \, \frac{B_{2j}}{(2j)!} \, \frac{1}{(n-2k)! \, (n-2j)! \, (2n-2k-2j+1)} } } \, .
\label{eq:Kuo3}
\end{eqnarray}
Now, put $(-1)^{n-1}$ in evidence and multiply both sides by $(2n)!$. This yields
\begin{eqnarray}
B_{2n} = (-1)^{n-1} \, \Biggl[ \frac{(2n)!}{4 \, {(n-1)!}^2 \, (2n-1)} -\frac{(2n)!}{(n-1)!} \, \sum_{k=0}^{\lfloor n/2 \rfloor}{\frac{B_{2k}}{(2k)!} \, \frac{1}{(n-2k)! \, (2n-2k)}} \nonumber \\
+ (2n)! \, \sum_{k=0}^{\lfloor n/2 \rfloor}{\sum_{j=0}^{\lfloor n/2 \rfloor}{\frac{B_{2k}}{(2k)!} \, \frac{B_{2j}}{(2j)!} \, \frac{1}{(n-2k)! \, (n-2j)! \, (2n-2k-2j+1)} } } \, ,
\label{eq:Kuo4}
\end{eqnarray}
which readily simplifies to
\begin{eqnarray}
B_{2n} = (-1)^{n-1} \, \Biggl[ \frac{n \, (2n-2)!}{2 \, (n-1)!^2} -\frac{(2n)!}{2\,(n-1)!}  \, \sum_{k=0}^{\lfloor n/2 \rfloor}{\frac{B_{2k}}{(2k)! \, (n-2k)! \, (n-k)}} \nonumber \\
 +(2n)! \, \sum_{k=0}^{\lfloor n/2 \rfloor}{\frac{B_{2k}}{(2k)! \, (n-2k)!} \, \sum_{j=0}^{\lfloor n/2 \rfloor}{\frac{B_{2j}}{(2j)! \, (n-2j)!} \: \frac{1}{2n-2k-2j+1} } } \Biggl] \, .
\end{eqnarray}
\begin{flushright} $\Box$ \end{flushright}
\end{prova}

In fact, the auxiliary terms $a_n$ and $b_n$ in our Theorem~\ref{teo:B2n} can be cast in a more suitable form for computational purposes. For $a_n$, one has $a_1 = \frac12$ and
\begin{equation}
a_n = \frac{n}{2\,(n-1)!} \, \frac{(2n-2)!}{(n-1)!} = \frac{n}{2\,(n-1)!} \, \prod_{m=2}^n{(2n -m)} \, ,
\label{eq:an1}
\end{equation}
for all $n>1$. 
For $b_n$, one has $b_1 = 1$ and
\begin{eqnarray}
b_n &=& \frac{(2n) \, (2n-1) \cdot \, \ldots \, \cdot n \, (n-1)!}{2\,(n-1)!}  = n^2 \, (2n-1) \cdot \, \ldots \, \cdot (n+1) \nonumber \\
&=& n^2 \, \prod_{m=1}^{n-1}{(2n -m)} \, ,
\end{eqnarray}
for all $n>1$.

Since $B_0 = 1$ and the recurrence relation in Theorem~\ref{teo:B2n} has only the four basic numeric operations ($+$, $-$, $\times$, and $\div$), it is straightforward to show, by induction on $n$, that every $B_{2n}$, $n \ge 1$, is a rational number, though this is a well-known characteristic of these numbers (see, e.g., Ref.~\cite{livrao}).



\end{document}